%% file: xg.tex
\documentclass{amsart}

\usepackage{amsmath}
\usepackage{graphicx,color}

\input{epsf.tex}

\theoremstyle{theorem}
\newtheorem{theorem}{Theorem}[section]   
\newtheorem{property}[theorem]{Property}
\theoremstyle{definition}
\newtheorem{proposition}[theorem]{Proposition}  
\newtheorem{definition}[theorem]{Definition}   
\newtheorem{rem}[theorem]{Remark}        
\numberwithin{equation}{section}     
\begin{document}
\title{A Homotopy Theory for Graphs } 

\author{E.~Babson}
\address{Department of Mathematics\\
  University of Washington\\Seattle, WA}
\email{babson@math.washington.edu}

\author{H.~Barcelo}
\address{Department of Mathematics and Statistics\\
  Arizona State University\\
  Tempe, Arizona 85287--1804} \email{barcelo@asu.edu}

\author{M.~De Longueville}
\address{Fachbereich Mathematik\\Freie Universit\"at Berlin\\
  Arnimallee 3--5, D-14195 Berlin, Germany}
\email{delong@math.fu-berlin.de}

\author{R.~Laubenbacher}
\address{Virginia Bioinformatics Institute\\Virginia 
Polytechnic Institute and State University\\Blacksburg, VA 24061}
\email{reinhard@vbi.vt.edu}

\date{}

\maketitle


\thispagestyle{empty}
\section{Introduction}

\begin{figure}[h]
  \begin{center}
    \includegraphics[width=0.25\textwidth]{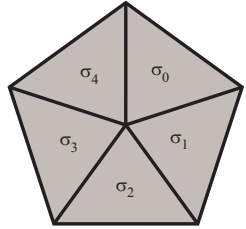}
    \caption{A $2$-dimensional complex $\Delta$ with nontrivial $A_1^1$.}
    \label{pentagon}
  \end{center}
\end{figure}

In the recent article \cite{BKLW} a new homotopy theory for graphs and
simplicial complexes was defined.  The motivation for the definition
came initially from a desire to find invariants for dynamic processes
that could be encoded via (combinatorial) simplicial complexes.  The
invariants should be topological in nature, but should at the same
time be sensitive to the combinatorics encoded in the complex, in
particular the level of connectivity of simplices (see \cite{KL}).
The construction is based on an approach proposed by R. Atkin
\cite{A1,A2}; hence the letter ``A.''  Namely, let $\Delta$ be a
simplicial complex of dimension $d$, let $0\leq q\leq d$ be an
integer, and let $\sigma_0\in\Delta$ be a simplex of dimension greater
than or equal to $q$.  One obtains a family of groups
$$
A_n^q(\Delta,\sigma_0),\quad n\geq 1,
$$ 
the $A$-groups of $\Delta$, based at $\sigma_0$.  These groups
differ from the classical homotopy groups of $\Delta$ in a 
significant way.  For instance, the group $A_1^1(\Delta,\sigma_0)$,
for the $2$-dimensional complex
$\Delta$ in Figure \ref{pentagon} is isomorphic to $\Bbb Z$, measuring
the presence of a ``connectivity'' hole in its center.  (See
the example on p. 101 of \cite{BKLW}.)  

\medskip
\par

The computation of these groups proceeds via the construction of a
graph, $\Gamma_q(\Delta)$, whose vertices represent simplices in
$\Delta$.  There is an edge between two simplices if they share a face
of dimension greater than or equal to $q$.  This construction
suggested a natural definition of the $A$-theory of graphs, which was
also developed in \cite{BKLW}.  Proposition 5.12 in that paper shows that
$A_1$ of the complex can be obtained as the fundamental group of the
space obtained by attaching $2$-cells into all $3$- and $4$-cycles of
$\Gamma_q(\Delta)$.

The goal of the present paper is to generalize this result.  Let
$\Gamma$ be a simple, undirected graph, with distinguished base vertex
$v_0$.  We will construct an infinite cell complex $X_{\Gamma}$
together with a homomorphism
$$
A_n(\Gamma,v_0)\longrightarrow \pi_n(X_{\Gamma},v_0).
$$\par
Moreover, we can show this homomorphism to be an isomorphism if a
(plausible) cubical analog of the simplicial approximation theorem
holds.

There are several reasons for this generalization.  One reason is the
desire for a homology theory associated to the $A$-theory of a graph.
A natural candidate is the singular homology of the space
$X_{\Gamma}$.  This will be explored in a future paper.

Another reason is a connection to the homotopy of the complements of
certain subspace arrangements.  While computing $A_1^{n-3}$ of the
order complex of the Boolean lattice $B_n$, it became clear that this
computation was equivalent to computing the fundamental group of the
complement of the $3$-equal arrangement \cite{BW}.  (This result for the
$k$-equal arrangement was proved independently by A. Bj\"orner \cite{Bj}.)
To generalize this connection to a wider class of subspace
arrangements a topological characterization of $A$-theory is needed.
\par

The content of the paper is as follows.  After a brief review of the
definition of $A$-theory, we construct the model space $X_{\Gamma}$,
followed by a proof of the main result (Theorem \ref{xg}).  The main
result refers to a yet unknown analog of a simplicial approximation
theorem in the cubical world (Property \ref{cat}), which we briefly
discuss in Section \ref{sec:cubicalstuff}. The last section introduces
the loop graph of a graph, and we prove that the $(n+1)$-st $A$-group
of the graph is isomorphic to the $n$-th $A$-group of the loop graph,
in analogy to a standard result about classical homotopy.

\section{$A$-theory of Graphs}

We first recall the definition given in Sect. 5 of \cite{BKLW}.

\begin{definition}{\rm
Let $\Gamma_1=(V_1,E_1),\ \Gamma_2=(V_2,E_2)$ be simple graphs, that is,
graphs without loops and multiple edges.
\begin{enumerate}
\item
The {\it Cartesian product} $\Gamma_1\times\Gamma_2$
is the graph with vertex set $V_1\times V_2$.  There is an
edge between $(u_1,u_2)$ and $(v_1,v_2)$ if either $u_1=v_1$
and $u_2v_2\in E_2$ or $u_2=v_2$ and $u_1v_1\in E_1$.
\item
A {\it graph homomorphism} $f:\Gamma_1\longrightarrow \Gamma_2$ is a
set map $V_1\longrightarrow V_2$ such that, if $uv\in E_1$,
then either $f(u)=f(v)$ or $f(u)f(v)\in E_2$.
\item
Let ${\bf I}_n$ be the graph with $n+1$ vertices labeled
$0,1,\ldots, n$, and edges $(i-1)i$ for $i=1,\ldots ,n$.
\item
Let $v_1\in\Gamma_1,v_2\in\Gamma_2$ be distinguished base vertices.
A {\it based} graph homomorphism $f:\Gamma_1\longrightarrow \Gamma_2$ is
a graph homomorphism such that $f(v_1)=v_2$.
\end{enumerate}
}
\end{definition}

Next we define homotopy of graph maps and homotopy equivalence of graphs.

\begin{definition}
{\rm

\begin{enumerate}
\item Let $f,g:(\Gamma_1,v_1)\longrightarrow (\Gamma_2,v_2)$ be based
  graph homomorphisms.  We call $f$ and $g$ {\it A--homotopic},
  denoted by $f\simeq_A g$, if there is an integer $n$ and a graph
  homomorphism
$$
\phi: \Gamma_1\times {\bf I}_n\longrightarrow \Gamma_2,
$$
such that $\phi(-,0)=f$, and $\phi(-,n)=g$, and such that
$\phi(v_1,i)=v_2$ for all $i$.
\end{enumerate}
}
\end{definition}

\begin{definition}
\begin{enumerate}
\item
Let
$$
{\bf I}^n_m={\bf I}_m\times \cdots \times{\bf I}_m
$$
be the $n$-fold Cartesian product of ${\bf I}_m$ for
some $m$.  We will call ${\bf I}^n_m$ an $n$-{\it cube} of
{\it height} $m$.  Its {\it distinguished base point} is
${\bf O}=(0,\ldots ,0)$.
\item
Define the {\it boundary} $\partial {\bf I}^n_m$
of a cube ${\bf I}^n_m$
of height $m$ to be the subgraph
of ${\bf I}^n_m$ containing all vertices with at least one
coordinate equal to $0$ or $m$.
\end{enumerate}
\end{definition}

It is easy to show (Lemma 5.4 of \cite{BKLW}) that any graph homomorphism
from ${\bf I}_m^n$ to $\Gamma$ can be extended to a graph homomorphism
from ${\bf I}_p^n$ to $\Gamma$ for any $p\geq m$.  Thus, by abuse of
notation we will sometimes omit the subscript $m$.

\begin{definition}{\rm
Let $A_n(\Gamma,v_0),\ n\geq 1$, be the set of homotopy
classes of graph homomorphisms 
$$
f:({\bf I}^n,\partial {\bf I}^n)\longrightarrow (\Gamma,v_0).
$$
For $n=0$, we define $A_0(\Gamma,v_0)$ to be the pointed
set of connected components of $\Gamma$, with distinguished
element the component containing $v_0$.  We will denote the
equivalence class of a homomorphism $f$ in $A_n(\Gamma,v_0)$ by $[f]$.

We can define a multiplication on the set
$A_n(\Gamma,v_0), n\geq 1$, as follows.  Given
elements $[f],\ [g]\in A_n(\Gamma,v_0)$, represented by
$$
f,g:({\bf I}^n_m,\partial {\bf I}^n_m)\longrightarrow (\Gamma,v_0),
$$
defined on a cube of height $m$,
we define $[f]*[g]\in A_n(\Gamma,v_0)$ as the homotopy class
of the map
$$
h:({\bf I}^{n}_{2m},\partial {\bf I}^{n}_{2m})\longrightarrow (\Gamma,v_0),
$$
defined on a cube of height $2m$ as follows.
$$
h(i_1,\ldots ,i_n)=\left\{\begin{array}{ll}
                     f(i_1,\ldots ,i_n) &
                          \mbox{if $i_j\leq m$ for all $j$},\\
                     g(i_1-m,\ldots ,i_n) &
                          \mbox{if $i_1>m$ and $i_j\leq m$ for $j>1$},\\
                     v_0                &\mbox{otherwise}.
                          \end{array}
               \right.
$$
}
\end{definition}

Alternatively, using Theorem 5.16 in \cite{BKLW}, one can describe the $A$-theory of graphs
using multidimensional ``grids'' of vertices as follows.  Let
$\Gamma$ be a graph with distinguished vertex $v_0$.  Let
$\mathcal A_n(\Gamma,v_0)$ be the set of functions
$$
\Bbb Z^n\longrightarrow V(\Gamma),
$$
from the lattice $\Bbb Z^n$ into the set of vertices of $\Gamma$
which take on the value $v_0$ almost everywhere, and for
which any two adjacent lattice points get mapped into either
the same or adjacent
vertices of $\Gamma$.  We define an
equivalence relation on this set as follows.  Two functions
$f$ and $g$ are equivalent, if there exists 
$$
h:\Bbb Z^{n+1}\longrightarrow V(\Gamma),
$$
in $\mathcal A_{n+1}(\Gamma,v_0)$ and integers $k$ and $l$, such that
\begin{align*}
h(i_1,\ldots ,i_n,k)&=f(i_1,\ldots ,i_n),\\
h(i_1,\ldots ,i_n,l)&=g(i_1,\ldots ,i_n) 
\end{align*}
for all $i_1,\ldots ,i_n\in \Bbb Z$.
For a definition of a group operation on the set of equivalence classes
see Prop. 3.5 of \cite{BKLW}.
Then it is straightforward to see that
$A_n(\Gamma,v_0)$ is isomorphic to the group of equivalence classes of elements
in $\mathcal A_n(\Gamma,v_0)$.  It will be useful to think of $A_n(\Gamma,v_0)$
in those terms.

\section{A cubical set setting for the $A$-theory of graphs}
We now define a cubical set $K_*(\Gamma)$ associated to the graph
$\Gamma$ (see \cite{Ma}). This gives the right setup in order to
obtain a close connection to the space $X_\Gamma$ which we define in
the next section. Let $I_{\infty}^n$ be the ``infinite'' discrete
$n$-cube, that is, the infinite lattice labeled by $\Bbb Z^n$.
\begin{definition}
  A graph homomorphism $f:I^n_\infty\rightarrow \Gamma$ {\em
    stabilizes in direction $(i,\varepsilon)$, $i=1,\ldots,n$,
    $\varepsilon\in\{\pm 1\}$} if there exists an $m_0$, s.t. for all
  $m\geq m_0$ 
$$f(a_1,\ldots,a_{i-1},\varepsilon m_0,a_{i+1},\ldots,a_n)=f(a_1,\ldots,a_{i-1},\varepsilon m,a_{i+1},\ldots,a_n).$$

\end{definition}
Let
$$
K_n(\Gamma)={\rm Hom}_s(I_{\infty}^n,\Gamma),
$$
the set of graph homomorphisms from the infinite $n$-cube to $\Gamma$
that eventually stabilize in each direction $(i,\varepsilon)$. 

For each ``face'' of $I_{\infty}^n$, i.e., for each choice of
$(i,\varepsilon)$, $i=1,\ldots,n$, $\varepsilon\in\{\pm 1\}$, we define {\em face
maps}
$$
\alpha_{i,\varepsilon}':K_n(\Gamma)\longrightarrow K_{n-1}(\Gamma),
$$
by 
$$\alpha_{i,\varepsilon}'(f)(a_1,\ldots,a_{n-1})=f(a_1,\ldots,a_{i-1},\varepsilon m_0,a_i,\ldots,a_{n-1}),$$
where $m_0$ is chosen large enough.
In other words $ \alpha_{i,\varepsilon}'(f)$ is the map in
$K_{n-1}(\Gamma)$ whose values are equal to the stable values of $f$
in direction $(i,\varepsilon)$.

{\em Degeneracy maps}
$$
\beta_i':K_{n-1}(\Gamma)\longrightarrow K_n(\Gamma),
$$
$i=1,\ldots ,n$, are defined as follows.  Given a map $f\in K_{n-1}(\Gamma)$,
extend it to a map on $I_{\infty}^n$ by 
$$
\beta_i'(f)(a_1,\ldots ,a_{n})=f(a_1,\ldots ,a_{i-1},a_{i+1},\ldots, a_n),
$$
for each $(a_1,\ldots ,a_{n})\in I_{\infty}^{n}$.
It is straightforward to check that in this way $K_*(\Gamma)$ is a cubical set.

We now imitate the definition of combinatorial homotopy of Kan complexes; see,
e.g., \cite[Ch. 1.3]{Ma}.

\begin{figure}[htbp]
  \centering
  \input{kan-equivalence.pstex_t}
  \caption{An illustration of a map $h$ in the definition of $\sim$.}
  \label{fig:kan-equivalence}
\end{figure}
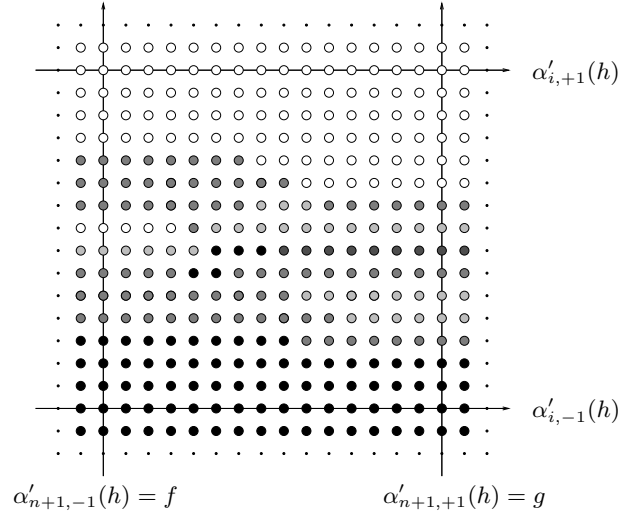

\begin{definition}
  We define a relation on $K_n(\Gamma)$, $n\geq 0$.  Let $f,g\in
  K_n(\Gamma)$. Then $f\sim g$ if there exists $h\in
  K_{n+1}(\Gamma)$ such that for all $i=1,\ldots,n$, $\varepsilon\in\{\pm1\}$:
\begin{enumerate}
\item
$\alpha_{i,\varepsilon}'(f)=\alpha_{i,\varepsilon}'(g)$,
\item
$\alpha_{i,\varepsilon}'(h)=\beta_{n}'\alpha_{i,\varepsilon}'(f)=
\beta_n'\alpha_{i,\varepsilon}'(g)$,
\item
 $\alpha_{n+1,-1}'(h)=f$ and
$\alpha_{n+1,1}'(h)=g$.
\end{enumerate}
For an illustration see Figure
  \ref{fig:kan-equivalence}. 
\end{definition}

\begin{proposition}
The relation defined above is an equivalence relation.\qed
\end{proposition}

\begin{definition}
Let $v_0\in \Gamma$ be a distinguished vertex.  
Let $B_*(\Gamma, v_0)\subset K_*(\Gamma)$ be the subset of all
maps that are equal to $v_0$ outside of a finite
region of $I_{\infty}^*$.   
\end{definition}

Observe that the equivalence relation $\sim$ restricts to an equivalence
relation on $B_*(\Gamma,v_0)$, also denoted by $\sim$. 

\begin{proposition}
There is a group structure on the 
set $B_n(\Gamma,v_0)/\sim$ for all $n\geq 1$, and, furthermore,
$$
(B_n(\Gamma,v_0)/\sim)\cong A_n(\Gamma, v_0).
$$
\end{proposition}\par
The proof is tedious, but straightforward.
For a definition of the group structure see Prop. 3.5 of \cite{BKLW}.
\section{Definition of $X_{\Gamma}$}

Let $\Gamma$ be a finite, simple (undirected) graph.  In this section
we define a cell complex $X_{\Gamma}$ associated to $\Gamma$.
This complex will be defined as the geometric realization of a 
certain cubical set $M_*(\Gamma)$.  Let $I^n_1$ be the discrete
$n$-cube.  Let
$$
M_n(\Gamma)={\rm Hom}(I_1^n,\Gamma),
$$
the set of all graph morphisms from $I_1^n$ to $\Gamma$.  We define
face and degeneracy maps as follows.

First note that $I_1^n$ has $2n$ faces $F_{i,\varepsilon}$, with $i=1,\ldots ,n$,
and $\varepsilon\in\{\pm 1\}$, corresponding to the two faces for each
coordinate.
For $i=1,\ldots ,n$, $\varepsilon\in\{\pm1\}$, let 
\begin{align*}
a_{i,\varepsilon}:I_1^{n-1}&\longrightarrow I_1^n\\
(x_1,\ldots,x_{n-1})&\longmapsto (x_1,\ldots,x_{i-1},{\textstyle\frac{\varepsilon+1}{2}},x_i,\ldots x_{n-1})
\end{align*}
be the graph map given by inclusion of $I_1^{n-1}$ as the
$(i,\varepsilon)$-face of $I_1^n$. For $i=1,\ldots ,n$ define
\begin{align*}
  b_i:I_1^n&\longrightarrow I_1^{n-1}\\
(x_1,\ldots,x_n)&\longmapsto (x_1,\ldots,x_{i-1},x_{i+1},\ldots x_n)
\end{align*}
to be the projection in direction $i$.

Now let
$$
\alpha_{i,\varepsilon}:M_n(\Gamma)\longrightarrow M_{n-1}(\Gamma)
$$
be the map induced by $a_{i,\varepsilon}$.  Likewise, define
$$
\beta_i:M_{n-1}(\Gamma)\longrightarrow M_n(\Gamma)
$$
to be the map induced by $b_i$.  In this way we obtain a cubical set
$M_*(\Gamma)$.

To each cubical set is associated a cell complex, namely its geometric
realization.  We recall the construction for $M_*(\Gamma)$.  Let $C^n$
be the geometric $n$-dimensional cube.  We can define functions
$a_{i,\varepsilon}$ and $b_i$ on $C^n$ in a fashion similar to above.
Define
the space
$$
|M_*(\Gamma)|=\biguplus_{n\geq 0}M_n(\Gamma)\times C^n/\sim ,
$$
where $\sim$ is the equivalence relation generated by the following
two types of equivalences:
\begin{align}
(\alpha_{i,\varepsilon}(f),x_{n-1})& \sim (f,a_{i,\varepsilon}
(x_{n-1})), \, f\in M_n(\Gamma),\, x_{n-1}\in C^{n-1}\\
(\beta_j(g),x_n)& \sim (g,b_j(x_n)),\, g\in M_{n-1}(\Gamma), \, x_n\in C^n.
\end{align}

We will denote the cell complex $|M_*(\Gamma)|$ by $X_{\Gamma}$.


\section{The main result}

We can now state the main result of the paper.

\begin{theorem}\label{xg}
There is a group homomorphism
$$
\phi:A_n(\Gamma, v_0)\longrightarrow \pi_n(X_{\Gamma},v_0),
$$
for all $n\geq 1$. If a cubical analog of the simplicial
approximation theorem such as \ref{cat} holds, then $\phi$ is an
isomorphism.
\end{theorem}

\begin{proof}
First we define $\phi$. Let $[f]\in A_n(\Gamma,v_0)\cong B_n(\Gamma,v_0)/\sim$.
Then a representative $f$ is a graph homomorphism
$$
f:I_{\infty}^n\longrightarrow \Gamma,
$$
whose value on vertices outside a finite region is equal to $v_0$,
say for vertices outside of a cube with side length $r$.  Our goal is
to define a continuous map
$$
\tilde f:C^n\longrightarrow X_{\Gamma},
$$
such that $\tilde f$ sends the boundary of $C^n$ to $v_0$.

Let $D^n$ be a cubical subdivision of $C^n$ into cubes of side
length $1/r$.  The 1-skeleton of $D^n$ can be identified with
$I_r^n$, which is contained in $I_{\infty}^n$.  And each
subcube of $I_r^n$ can be identified with $I_1^n$.  Hence, $f$
restricts to a graph homomorphism on each cube in the 1-skeleton of
$D^n$, that is, a graph homomorphism
$$
\hat f:I_1^n\longrightarrow \Gamma.
$$
Thus, $\hat f\in {\rm Hom}(I_1^n,\Gamma)$.  Now define
$\tilde f$ on each subcube of $D^n$ by
$$
\tilde f(x)=[(\hat f,x)]\in X_{\Gamma}=(\biguplus_n{\rm Hom}(I_1^n,\Gamma)\times C^n)/\sim.
$$
The equivalence relation $\sim$ guarantees that $\tilde f$ is well-defined on overlapping
faces.  Therefore, our definition extends to give a map
$$
\tilde f:D^n\longrightarrow X_{\Gamma}.
$$
So define 
$$
\phi([f])=[\tilde f].
$$

We need to show that $\phi$ is well-defined.  Let $f\sim g$ be two maps in
$B_n(\Gamma, v_0)$.  Then there exists a homotopy $h\in B_{n+1}(\Gamma ,v_0)$
such that $\alpha_{n+1,-1}'(h)=f$ and $\alpha_{n+1,1}'(h)=g$.  We claim that
$\phi(h)$ gives a homotopy between $\phi(f)$ and $\phi(g)$.  
From the definition of $\phi$ it is easy to see that
$$
\phi((\alpha_{i,\varepsilon}'(h))(y)=[(\alpha_{i,\varepsilon}(\tilde h),y)],
$$
for all $i,\varepsilon$.  Therefore, the restriction of 
$$
\phi(h):D^{n+1}\longrightarrow X_{\Gamma}
$$
to the $(n+1,-1)$-face is equal to the map from
$D^n$ to $X_{\Gamma}$,
sending $x$ to $[(\alpha_{n+1,-1}'(h),x)]$, which is equal to $\phi(f)$;
similarly for $\phi(g)$.
It now follows that $\phi(h)$ is a homotopy
between $\phi(f)$ and $\phi(g)$.  This shows that $\phi$ is well-defined.

Now we show that $\phi$ is a group homomorphism.  Recall \cite[p.~111]{BKLW} that the multiplication in $A_n(\Gamma,v_0)$ is given by
juxtaposing ``grids.''  This carries over directly to
$B_n(\Gamma,v_0)/\sim$.  On the other hand, the multiplication in
$\pi_n(X_{\Gamma},v_0)$ is given by using the comultiplication on
$(C^n,\partial C^n)$.  It is then straightforward to check that $\phi$
preserves multiplication.  

From here on we assume that Property \ref{cat} holds. Under this
assumption we show that $\phi$ is onto.  We first show that every
element in $\pi_n(X_{\Gamma},v_0)$ contains a cubical representative.
Let $[f]\in\pi_n(X_{\Gamma}, v_0)$.  Then $f:C^n\longrightarrow
X_{\Gamma}$ sends the boundary of $C^n$ to the base point $v_0$.
Trivially then, the restriction of $f$ to the boundary is a cubical
map.  By Property \ref{cat} $f$ is homotopic to a cubical map on a
cubical subdivision $D^n$ of $C^n$, and agrees with $f$ on the
boundary.  That is, $[f]$ contains a cubical representative.  So we
may assume that $f$ is cubical on $D^n$.

Consider the restriction of $f$ to the $1$-skeleton of $D^n$.  
It induces in the obvious way a graph map $g:I_{\infty}^n\longrightarrow \Gamma$,
that is, an element $[g]\in B_n(\Gamma,v_0)/\sim $.  
We claim that $\phi(g)=[f]$, that is, $\tilde g\sim f$.  We use
induction on $n$.  If $n=1$, then we are done, since any two maps on the
unit interval that agree on the end points are homotopic.
Changing $f$ up to homotopy we may assume that $f$ and $\tilde g$ are equal on
the $1$-skeleton.
 
Now let $n>1$.  Note that 
$$
f:D^n\longrightarrow X_{\Gamma}=(\biguplus_{n\geq 0}{\rm Hom}(I^n_1,\Gamma))\times C^n/\sim 
$$
is cubical, so each $n$-cube $C^n$ in the cubical subdivision $D^n$ is
sent to an $n$-cube in $X_{\Gamma}$.  The particular $n$-cube it is mapped to
is determined by the image of the map on the $1$-skeleton, since the map
is cubical.  This in turn determines an element in ${\rm Hom}(I_1^n,\Gamma)$,
serving as the label of the image cube.  Hence, $f$ and $\tilde g$ map each
$n$-cube of the subdivision $D^n$ to the same $n$-cube in $X_{\Gamma}$.  
By induction we may assume that $f$ and $\tilde g$ are equal on the boundary
of each $n$-cube.  But observe that any two maps into $C^n$ that agree on
the boundary are homotopic, via a homotopy that leaves the boundary fixed.
This shows $f$ and $\tilde g$ are homotopic on each $n$-cube of the 
cubical subdivision $D^n$.  Pasting these homotopies together along the
boundaries, we obtain a homotopy between $f$ and $\tilde g$, so that 
$[f]=[\tilde g]$. 

To show that $\phi$ is one-to-one under the assumption of Property
\ref{cat}, suppose that $f,g\in B_n(\Gamma,v_0)/\sim$ such that
$\phi(f)=\phi(g)\in \pi_n(\Gamma,v_0)$.  Then there exists a homotopy
$h:C^{n+1}\longrightarrow X_{\Gamma}$ such that the restrictions of
$h$ to the $(n+1)$-directional faces are $\phi(f)$ and $\phi(g)$,
respectively.  As above, we may assume that $h$ is cubical on a
subdivision $D^{n+1}$ of $C^{n+1}$, providing a homotopy between
cubical approximations of $\phi(f)$ and $\phi(g)$ on a subdivision
$D^n$ of $C^n$.  Now observe that the restriction of $h$ to the
$1$-skeleton of $D^{n+1}$ induces a graph homomorphism
$h':I^{n+1}_1\longrightarrow \Gamma$ in $B_{n+1}(\Gamma,v_0)$, whose
restrictions to the $(n+1)$-directional faces are refinements of $f$
and $g$, respectively.  But these refinements are equivalent to $f$
and $g$, respectively.  Thus, $[f]=[g]\in B_n(\Gamma,v_0)/\sim$. 
\end{proof}

\section{Cubical Complexes}
\label{sec:cubicalstuff}

The following plausible property is a special case of a general
cubical approximation theorem. We have not found it in the literature and
have not been able to prove it yet.

\begin{property}\label{cat}
Let $X$ be a cubical set, and let $f:C^n\longrightarrow |X|$
be a continuous map from the $n$-cube to the geometric
realization of $X$, such that the restriction of $f$ to the
boundary of $C^n$ is cubical.  Then there exists a cubical
subdivision $D^n$ of $C^n$ and a cubical map $f':D^n\longrightarrow |X|$
which is homotopic to $f$ and the restrictions of $f$ and $f'$ to
the boundary of $D^n$ are equal.  
\end{property}

\section{Path- and loop graph of a graph}
In topology the computation of the homotopy group $\pi_{n+1}(X)$ of a
space $X$ can be reduced to the computation of $\pi_{n}(\Omega X)$,
the $n$-th homotopy group of the loop space $\Omega X$ of $X$. Here we
want to introduce the path graph $PG$ and the loop graph $\Omega G$ of
a graph $G$ such that naturally  $A_n(\Omega G)\cong A_{n+1}(G) $.
\begin{definition}
  Let $G$ be a graph with base vertex $\ast$.  Define the {\em path graph} $PG=(V_{PG},E_{PG})$
  to be the graph on the vertex set
  \begin{align*}
    V_{PG}=\{\varphi:I_m\rightarrow G:m\in \mathbb N,\,\varphi\text{ a graph map with } \varphi(0)=\ast\}.
  \end{align*}
  The edge set $E_{PG}$ is given as follows. Consider two vertices
  $\varphi_0:I_m\rightarrow G$ and $\varphi_1:I_{m'}\rightarrow G$.
  Assuming $m\leq m'$ extend $\varphi_0$ to a map
  $\varphi_0':I_{m'}\rightarrow G$ by repeating the last vertex
  $\varphi_0(m)$ at the end:
\begin{align*}
\varphi_0'(y)=
\begin{cases}
  \varphi_0(y), &\text{ if } y\leq m,\\
  \varphi_0(m), &\text{ otherwise.}
\end{cases} 
\end{align*}
Define $\{\varphi_0,\varphi_1\}$ to be an edge if there exists a graph map 
$\Phi:I_{m'}\times I_1\rightarrow G$ such that $\Phi(\bullet,0)=\varphi_0'$ and $\Phi(\bullet,1)=\varphi_1$.
\end{definition}
There is graph map $p:PG\rightarrow G$ given by
$p(\varphi)=\varphi(m)$ for a vertex $\varphi:I_m\rightarrow G$ of
$PG$.
\begin{definition}
  For a graph $G$ define the {\em loop graph} $\Omega G$ of $G$ to be the
  induced subgraph of $PG$ on the vertex set $p^{-1}(\ast)$. We define
  the base vertex of $\Omega G$ to be the vertex $\varphi_0:I_0\rightarrow
  G$, i.e., the map that sends the single vertex of $I_0$ to $\ast$ in
  $G$. To avoid too much notation we will denote this map by $\ast$ as well.
\end{definition}
Note that for a graph map $\psi :(G,\ast)\rightarrow(H,\ast)$ there is
an induced map $\Omega \psi :(\Omega G,\ast)\rightarrow (\Omega
H,\ast)$ defined by $\Omega \psi (\varphi)(y)=\psi (\varphi(y))$ where
$\varphi:I_m\rightarrow G$ and $y$ is a vertex of $I_m$.
\begin{rem}
\label{constantloop}
Consider the constant loop $\varphi_m:I_m\rightarrow G$ in $\Omega G$,
i.e., $\varphi_m(x)=\ast\in G$ for all vertices $x$ of $I_m$. If a
loop $\varphi:I_m\rightarrow G$ is connected to $\varphi_m$ via an
edge, then it is also connected to $\varphi_0=\ast$ via an edge.
\end{rem}
Analogously to classical topology we have the following.
\begin{proposition}
  There is a natural isomorphism $A_n(\Omega
  G)\xrightarrow{\cong}{}A_{n+1}(G)$ for $n\geq 1$. Furthermore, there
  is a bijection $A_0(\Omega G)\xrightarrow{\cong}{}A_{1}(G)$.
\end{proposition}
\begin{proof}
  {\bf The case $n\geq 1$.}  Let $[f]\in A_n(\Omega G)$, i.e., $f$ is
  a graph map $f:(I_m^n,\partial I_m^n)\rightarrow (\Omega G,\ast)$.
  For $x$ a vertex of $I_m^n$ there is an $m_f(x)$ such that $f(x)$ is
  a graph map $f(x):(I_{m_f(x)},\partial I_{m_f(x)})\rightarrow
  (G,*)$. Let $m'=\max_{x}\{m_f(x),m\}$. We want to define a graph map
  $\alpha(f):(I_{m'}^{n+1},\partial I_{m'}^{n+1})\rightarrow
  (G,\ast)$. For that reason write $I_{m'}^{n+1}=I_{m'}^n\times
  I_{m'}$ and let $(x,y)$ be a vertex of $I_{m'}^n\times I_{m'}$. Now
  let
  \begin{align*}
    \alpha(f)(x,y)=
    \begin{cases}
      f(x)(y), &\text{ if $x$ is a vertex of $I^n_m\subset I^n_{m'}$ and } y\leq m_f(x),\\
      \ast, &\text{ otherwise.}
\end{cases}
  \end{align*}
  The construction is shown in Figure \ref{fig:loop-iso}, where $n=1$,
  $m=10$, and $m'=12$. The vertical line is $I^n_m$, the horizontal
  lines indicate the paths $f(x)$, the whole square indicates $\alpha
  (f)$.
  \begin{figure}[htbp]
    \begin{center}
      \input{loop-iso.pstex_t}
      \caption{The maps $f$  and $\alpha(f)$.}
      \label{fig:loop-iso}
    \end{center}
  \end{figure}
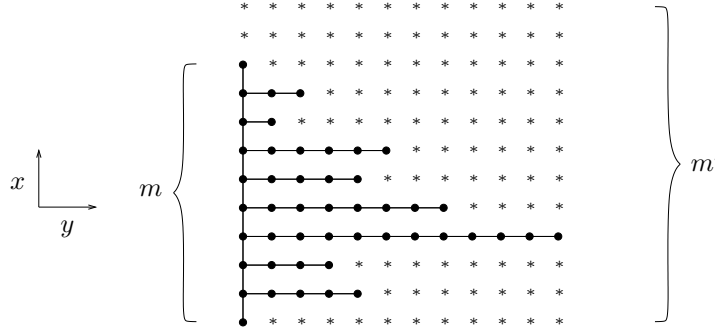\par
  We claim that the map $[f]\rightarrow [\alpha(f)]$ is well defined
  and the desired natural isomorphism.\par
  
  \noindent{\bf Well definedness:} First of all it is easy to check that
  $\alpha(f)$ is a graph map $\alpha(f):(I_{m'}^{n+1},\partial
  I_{m'}^{n+1})\rightarrow (G,\ast)$. Now let $[f]=[g]\in A_n(\Omega
  G)$, i.e., there exists an A--homotopy $H:I^n_m\times I_l\rightarrow
  \Omega G$ between $f$ and $g$. Now let
  $m'=\max_{x,x'}\{m_f(x),m_g(x'),m\}$ and define $\bar
  H:I^n_{m'}\times I_{m'}\times I_l\rightarrow G$ by
  \begin{align*}
    \bar H(x,y,t)=
    \begin{cases}
      H(x,t)(y), &\text{ if $x$ is a vertex of $I^n_m$ and $y\leq m_{H(x,t)}$,}\\
\ast, &\text{ otherwise.}
    \end{cases}
  \end{align*}
  Then $\bar H$ is a graph map and an A--homotopy between (possibly
  extended to a larger cube) $\alpha(f)$ and $\alpha(g)$.\par
  
  \noindent{\bf Homomorphism:} Is straightforward; similar techniques play a
  role that are needed to show that $A_{n}(G)$ is a group for $n\geq
  1$.\par

  \noindent{\bf Surjectivity:} For $[h]\in A_{n+1}(G)$, say
  $h:(I^{n+1}_m,\partial I^{n+1}_m)\rightarrow (G,\ast)$, consider the
  map $f$ defined by $f(x)(y)=h(x,y)$ for $x$ a vertex of $I^n_m, y$ a
  vertex of $I_m$. This map is not quite what we want since it is a
  map $f:(I^n_m,\partial I^n_m)\rightarrow (\Omega G, \varphi_m)$,
  where $\varphi_m$ is the constant loop $I_m\rightarrow G$ as in
  Remark \ref{constantloop}.  Now define $f':(I^n_{m},\partial
  I^n_{m})\rightarrow (\Omega G,\ast)$ by $f'(x)=\ast\in \Omega G$ for
  $x$ a vertex of $\partial I^n_{m}$ and $f'(x)=f(x)$ for $x$ a vertex
  of $I^n_{m}\setminus\partial I^n_{m}$.  Thanks to Remark
  \ref{constantloop}, $f'$ is a well defined graph map and clearly
  $\alpha(f')=h$.
  \par
  
  \noindent{\bf Injectivity:} Consider $f:(I^n_m,\partial I^n_m)\rightarrow
  (\Omega G,\ast)$ and $g:(I^n_{m'},\partial I^n_{m'})\rightarrow
  (\Omega G,\ast)$ such that $[\alpha(f)]=[\alpha(g)]$, i.e., there is
  an A--homotopy $H:I^n_{m''}\times I_l\rightarrow G$ between
  (possibly extended to a larger cube) $\alpha(f)$ and $\alpha(g)$,
  where $m''=\max_{x,x'}\{m_f(x),m_g(x'),m,m'\}$.  Define $\bar
  H:I^n_{m''}\times I_l\rightarrow\Omega G$ by $\bar
  H(x,t)(y)=H(x,y,t)$. Then $\bar H(x,t):I_{m''}\rightarrow \Omega G$
  for all $x$ and $t$. Furthermore $\bar H(x,t)=\varphi_m$ for $x$ a
  vertex of $\partial I^n_{m''}$. As before we replace $\bar H$ by
  $\bar H'$ by changing it only on the boundary and by replacing $\alpha(f)$ by $f$ and $\alpha(g)$ by $g$.
  \begin{align*}
    \bar H'(x,t)=
    \begin{cases}
      \bar H(x,t), &\text{ if } x\text{ a vertex of } I^n_{m''}\setminus\partial I^n_{m''} \text{ and }t\not=0,m,\\
      f(x), &\text{ if } t=0 \text{ and } x \text{ a vertex of } I^n_m\subset I^n_{m''},\\
      g(x), &\text{ if } t=m \text{ and } x \text{ a vertex of } I^n_{m'}\subset I^n_{m''},\\
      \varphi_0, &\text{ otherwise.}
    \end{cases}
  \end{align*}
  Then by Remark \ref{constantloop} $\bar H'$ is a graph map and it
  yields an A--homotopy between (possibly extended to a larger cube)
  $f$ and $g$.\par

\noindent{\bf Naturality:} Let $\psi :(G,\ast_G)\rightarrow (H,\ast_H)$ be a graph map
and $f:(I_m^n,\partial I_m^n)\rightarrow (\Omega G,\ast)$. Then for a
vertex $x$ of $I_m^n$ we obtain
  \begin{align*}
  \psi _\#(\alpha_G(f))(x,y)&=
  \begin{cases}
    \psi (f(x)(y)), &\text{ if }y\leq m_f(x),\\
      \psi (\ast_G), &\text{ otherwise.}    
  \end{cases}\\
&=\begin{cases}
    \Omega \psi(f)(x)(y), &\text{ if }y\leq m_{\Omega \psi (f)}(x),\\
      \ast_H, &\text{ otherwise.}    
  \end{cases}\\
&=\alpha_H((\Omega \psi )_\#(f)).
\end{align*}\par
\noindent {\bf The remaining case $n=0$:} 
Consider an element $[\varphi]$ of $A_0(\Omega G)$, i.e., a connected
component of $\Omega G$ represented by a loop $\varphi:I_m\rightarrow
G$. This loop defines an element $[\varphi]$ (this time a homotopy
class) of $A_1(G)$. Well definedness and bijectivity of this
assignment is immediate.
\end{proof}
\subsection*{Acknowledgements}
The authors thank Rick Jardine and Vic Reiner for several helpful conversations.

\end{document}

%% file: kan-equivalence.pstex_t
\begin{picture}(0,0)%
\includegraphics{kan-equivalence.pstex}%
\end{picture}%
\setlength{\unitlength}{1243sp}%
\begingroup\makeatletter\ifx\SetFigFont\undefined
\def\x#1#2#3#4#5#6#7\relax{\def\x{#1#2#3#4#5#6}}%
\expandafter\x\fmtname xxxxxx\relax \def\y{splain}%
\ifx\x\y   
\gdef\SetFigFont#1#2#3{%
  \ifnum #1<17\tiny\else \ifnum #1<20\small\else
  \ifnum #1<24\normalsize\else \ifnum #1<29\large\else
  \ifnum #1<34\Large\else \ifnum #1<41\LARGE\else
     \huge\fi\fi\fi\fi\fi\fi
  \csname #3\endcsname}%
\else
\gdef\SetFigFont#1#2#3{\begingroup
  \count@#1\relax \ifnum 25<\count@\count@25\fi
  \def\x{\endgroup\@setsize\SetFigFont{#2pt}}%
  \expandafter\x
    \csname \romannumeral\the\count@ pt\expandafter\endcsname
    \csname @\romannumeral\the\count@ pt\endcsname
  \csname #3\endcsname}%
\fi
\fi\endgroup
\begin{picture}(10350,10132)(901,-9281)
\put(11251,-7441){\makebox(0,0)[lb]{\smash{\SetFigFont{9}{10.8}{rm}{\color[rgb]{0,0,0}$\alpha_{i,-1}'(h)$}%
}}}
\put(11251,-691){\makebox(0,0)[lb]{\smash{\SetFigFont{9}{10.8}{rm}{\color[rgb]{0,0,0}$\alpha_{i,+1}'(h)$}%
}}}
\put(901,-9151){\makebox(0,0)[lb]{\smash{\SetFigFont{9}{10.8}{rm}{\color[rgb]{0,0,0}$\alpha_{n+1,-1}'(h)=f$}%
}}}
\put(8281,-9151){\makebox(0,0)[lb]{\smash{\SetFigFont{9}{10.8}{rm}{\color[rgb]{0,0,0}$\alpha_{n+1,+1}'(h)=g$}%
}}}
\end{picture}

%% file: loop-iso.pstex_t
\begin{picture}(0,0)%
\includegraphics{loop-iso.pstex}%
\end{picture}%
\setlength{\unitlength}{1184sp}%
\begingroup\makeatletter\ifx\SetFigFont\undefined
\def\x#1#2#3#4#5#6#7\relax{\def\x{#1#2#3#4#5#6}}%
\expandafter\x\fmtname xxxxxx\relax \def\y{splain}%
\ifx\x\y   
\gdef\SetFigFont#1#2#3{%
  \ifnum #1<17\tiny\else \ifnum #1<20\small\else
  \ifnum #1<24\normalsize\else \ifnum #1<29\large\else
  \ifnum #1<34\Large\else \ifnum #1<41\LARGE\else
     \huge\fi\fi\fi\fi\fi\fi
  \csname #3\endcsname}%
\else
\gdef\SetFigFont#1#2#3{\begingroup
  \count@#1\relax \ifnum 25<\count@\count@25\fi
  \def\x{\endgroup\@setsize\SetFigFont{#2pt}}%
  \expandafter\x
    \csname \romannumeral\the\count@ pt\expandafter\endcsname
    \csname @\romannumeral\the\count@ pt\endcsname
  \csname #3\endcsname}%
\fi
\fi\endgroup
\begin{picture}(14250,6864)(-1799,-5908)
\put(3601,-5836){\makebox(0,0)[lb]{\smash{\SetFigFont{6}{7.2}{rm}{\color[rgb]{0,0,0}$\ast$}%
}}}
\put(4201,-5836){\makebox(0,0)[lb]{\smash{\SetFigFont{6}{7.2}{rm}{\color[rgb]{0,0,0}$\ast$}%
}}}
\put(4801,-5836){\makebox(0,0)[lb]{\smash{\SetFigFont{6}{7.2}{rm}{\color[rgb]{0,0,0}$\ast$}%
}}}
\put(5401,-5836){\makebox(0,0)[lb]{\smash{\SetFigFont{6}{7.2}{rm}{\color[rgb]{0,0,0}$\ast$}%
}}}
\put(6001,-5836){\makebox(0,0)[lb]{\smash{\SetFigFont{6}{7.2}{rm}{\color[rgb]{0,0,0}$\ast$}%
}}}
\put(6601,-5836){\makebox(0,0)[lb]{\smash{\SetFigFont{6}{7.2}{rm}{\color[rgb]{0,0,0}$\ast$}%
}}}
\put(7201,-5836){\makebox(0,0)[lb]{\smash{\SetFigFont{6}{7.2}{rm}{\color[rgb]{0,0,0}$\ast$}%
}}}
\put(7801,-5836){\makebox(0,0)[lb]{\smash{\SetFigFont{6}{7.2}{rm}{\color[rgb]{0,0,0}$\ast$}%
}}}
\put(8401,-5836){\makebox(0,0)[lb]{\smash{\SetFigFont{6}{7.2}{rm}{\color[rgb]{0,0,0}$\ast$}%
}}}
\put(9001,-5836){\makebox(0,0)[lb]{\smash{\SetFigFont{6}{7.2}{rm}{\color[rgb]{0,0,0}$\ast$}%
}}}
\put(9601,-5836){\makebox(0,0)[lb]{\smash{\SetFigFont{6}{7.2}{rm}{\color[rgb]{0,0,0}$\ast$}%
}}}
\put(6601,-5236){\makebox(0,0)[lb]{\smash{\SetFigFont{6}{7.2}{rm}{\color[rgb]{0,0,0}$\ast$}%
}}}
\put(7201,-5236){\makebox(0,0)[lb]{\smash{\SetFigFont{6}{7.2}{rm}{\color[rgb]{0,0,0}$\ast$}%
}}}
\put(7801,-5236){\makebox(0,0)[lb]{\smash{\SetFigFont{6}{7.2}{rm}{\color[rgb]{0,0,0}$\ast$}%
}}}
\put(8401,-5236){\makebox(0,0)[lb]{\smash{\SetFigFont{6}{7.2}{rm}{\color[rgb]{0,0,0}$\ast$}%
}}}
\put(9001,-5236){\makebox(0,0)[lb]{\smash{\SetFigFont{6}{7.2}{rm}{\color[rgb]{0,0,0}$\ast$}%
}}}
\put(9601,-5236){\makebox(0,0)[lb]{\smash{\SetFigFont{6}{7.2}{rm}{\color[rgb]{0,0,0}$\ast$}%
}}}
\put(6001,-5236){\makebox(0,0)[lb]{\smash{\SetFigFont{6}{7.2}{rm}{\color[rgb]{0,0,0}$\ast$}%
}}}
\put(6601,-2836){\makebox(0,0)[lb]{\smash{\SetFigFont{6}{7.2}{rm}{\color[rgb]{0,0,0}$\ast$}%
}}}
\put(7201,-2836){\makebox(0,0)[lb]{\smash{\SetFigFont{6}{7.2}{rm}{\color[rgb]{0,0,0}$\ast$}%
}}}
\put(7801,-2836){\makebox(0,0)[lb]{\smash{\SetFigFont{6}{7.2}{rm}{\color[rgb]{0,0,0}$\ast$}%
}}}
\put(8401,-2836){\makebox(0,0)[lb]{\smash{\SetFigFont{6}{7.2}{rm}{\color[rgb]{0,0,0}$\ast$}%
}}}
\put(9001,-2836){\makebox(0,0)[lb]{\smash{\SetFigFont{6}{7.2}{rm}{\color[rgb]{0,0,0}$\ast$}%
}}}
\put(9601,-2836){\makebox(0,0)[lb]{\smash{\SetFigFont{6}{7.2}{rm}{\color[rgb]{0,0,0}$\ast$}%
}}}
\put(6001,-2836){\makebox(0,0)[lb]{\smash{\SetFigFont{6}{7.2}{rm}{\color[rgb]{0,0,0}$\ast$}%
}}}
\put(3601,-436){\makebox(0,0)[lb]{\smash{\SetFigFont{6}{7.2}{rm}{\color[rgb]{0,0,0}$\ast$}%
}}}
\put(4201,-436){\makebox(0,0)[lb]{\smash{\SetFigFont{6}{7.2}{rm}{\color[rgb]{0,0,0}$\ast$}%
}}}
\put(4801,-436){\makebox(0,0)[lb]{\smash{\SetFigFont{6}{7.2}{rm}{\color[rgb]{0,0,0}$\ast$}%
}}}
\put(5401,-436){\makebox(0,0)[lb]{\smash{\SetFigFont{6}{7.2}{rm}{\color[rgb]{0,0,0}$\ast$}%
}}}
\put(6001,-436){\makebox(0,0)[lb]{\smash{\SetFigFont{6}{7.2}{rm}{\color[rgb]{0,0,0}$\ast$}%
}}}
\put(6601,-436){\makebox(0,0)[lb]{\smash{\SetFigFont{6}{7.2}{rm}{\color[rgb]{0,0,0}$\ast$}%
}}}
\put(7201,-436){\makebox(0,0)[lb]{\smash{\SetFigFont{6}{7.2}{rm}{\color[rgb]{0,0,0}$\ast$}%
}}}
\put(7801,-436){\makebox(0,0)[lb]{\smash{\SetFigFont{6}{7.2}{rm}{\color[rgb]{0,0,0}$\ast$}%
}}}
\put(8401,-436){\makebox(0,0)[lb]{\smash{\SetFigFont{6}{7.2}{rm}{\color[rgb]{0,0,0}$\ast$}%
}}}
\put(9001,-436){\makebox(0,0)[lb]{\smash{\SetFigFont{6}{7.2}{rm}{\color[rgb]{0,0,0}$\ast$}%
}}}
\put(9601,-436){\makebox(0,0)[lb]{\smash{\SetFigFont{6}{7.2}{rm}{\color[rgb]{0,0,0}$\ast$}%
}}}
\put(3001,164){\makebox(0,0)[lb]{\smash{\SetFigFont{6}{7.2}{rm}{\color[rgb]{0,0,0}$\ast$}%
}}}
\put(3601,164){\makebox(0,0)[lb]{\smash{\SetFigFont{6}{7.2}{rm}{\color[rgb]{0,0,0}$\ast$}%
}}}
\put(4201,164){\makebox(0,0)[lb]{\smash{\SetFigFont{6}{7.2}{rm}{\color[rgb]{0,0,0}$\ast$}%
}}}
\put(4801,164){\makebox(0,0)[lb]{\smash{\SetFigFont{6}{7.2}{rm}{\color[rgb]{0,0,0}$\ast$}%
}}}
\put(5401,164){\makebox(0,0)[lb]{\smash{\SetFigFont{6}{7.2}{rm}{\color[rgb]{0,0,0}$\ast$}%
}}}
\put(6001,164){\makebox(0,0)[lb]{\smash{\SetFigFont{6}{7.2}{rm}{\color[rgb]{0,0,0}$\ast$}%
}}}
\put(6601,164){\makebox(0,0)[lb]{\smash{\SetFigFont{6}{7.2}{rm}{\color[rgb]{0,0,0}$\ast$}%
}}}
\put(7201,164){\makebox(0,0)[lb]{\smash{\SetFigFont{6}{7.2}{rm}{\color[rgb]{0,0,0}$\ast$}%
}}}
\put(7801,164){\makebox(0,0)[lb]{\smash{\SetFigFont{6}{7.2}{rm}{\color[rgb]{0,0,0}$\ast$}%
}}}
\put(8401,164){\makebox(0,0)[lb]{\smash{\SetFigFont{6}{7.2}{rm}{\color[rgb]{0,0,0}$\ast$}%
}}}
\put(9001,164){\makebox(0,0)[lb]{\smash{\SetFigFont{6}{7.2}{rm}{\color[rgb]{0,0,0}$\ast$}%
}}}
\put(3001,764){\makebox(0,0)[lb]{\smash{\SetFigFont{6}{7.2}{rm}{\color[rgb]{0,0,0}$\ast$}%
}}}
\put(3601,764){\makebox(0,0)[lb]{\smash{\SetFigFont{6}{7.2}{rm}{\color[rgb]{0,0,0}$\ast$}%
}}}
\put(4201,764){\makebox(0,0)[lb]{\smash{\SetFigFont{6}{7.2}{rm}{\color[rgb]{0,0,0}$\ast$}%
}}}
\put(4801,764){\makebox(0,0)[lb]{\smash{\SetFigFont{6}{7.2}{rm}{\color[rgb]{0,0,0}$\ast$}%
}}}
\put(5401,764){\makebox(0,0)[lb]{\smash{\SetFigFont{6}{7.2}{rm}{\color[rgb]{0,0,0}$\ast$}%
}}}
\put(6001,764){\makebox(0,0)[lb]{\smash{\SetFigFont{6}{7.2}{rm}{\color[rgb]{0,0,0}$\ast$}%
}}}
\put(6601,764){\makebox(0,0)[lb]{\smash{\SetFigFont{6}{7.2}{rm}{\color[rgb]{0,0,0}$\ast$}%
}}}
\put(7201,764){\makebox(0,0)[lb]{\smash{\SetFigFont{6}{7.2}{rm}{\color[rgb]{0,0,0}$\ast$}%
}}}
\put(7801,764){\makebox(0,0)[lb]{\smash{\SetFigFont{6}{7.2}{rm}{\color[rgb]{0,0,0}$\ast$}%
}}}
\put(8401,764){\makebox(0,0)[lb]{\smash{\SetFigFont{6}{7.2}{rm}{\color[rgb]{0,0,0}$\ast$}%
}}}
\put(9001,764){\makebox(0,0)[lb]{\smash{\SetFigFont{6}{7.2}{rm}{\color[rgb]{0,0,0}$\ast$}%
}}}
\put(4201,-1636){\makebox(0,0)[lb]{\smash{\SetFigFont{6}{7.2}{rm}{\color[rgb]{0,0,0}$\ast$}%
}}}
\put(4801,-1636){\makebox(0,0)[lb]{\smash{\SetFigFont{6}{7.2}{rm}{\color[rgb]{0,0,0}$\ast$}%
}}}
\put(5401,-1636){\makebox(0,0)[lb]{\smash{\SetFigFont{6}{7.2}{rm}{\color[rgb]{0,0,0}$\ast$}%
}}}
\put(6001,-1636){\makebox(0,0)[lb]{\smash{\SetFigFont{6}{7.2}{rm}{\color[rgb]{0,0,0}$\ast$}%
}}}
\put(6601,-1636){\makebox(0,0)[lb]{\smash{\SetFigFont{6}{7.2}{rm}{\color[rgb]{0,0,0}$\ast$}%
}}}
\put(7201,-1636){\makebox(0,0)[lb]{\smash{\SetFigFont{6}{7.2}{rm}{\color[rgb]{0,0,0}$\ast$}%
}}}
\put(7801,-1636){\makebox(0,0)[lb]{\smash{\SetFigFont{6}{7.2}{rm}{\color[rgb]{0,0,0}$\ast$}%
}}}
\put(8401,-1636){\makebox(0,0)[lb]{\smash{\SetFigFont{6}{7.2}{rm}{\color[rgb]{0,0,0}$\ast$}%
}}}
\put(9001,-1636){\makebox(0,0)[lb]{\smash{\SetFigFont{6}{7.2}{rm}{\color[rgb]{0,0,0}$\ast$}%
}}}
\put(9601,-1636){\makebox(0,0)[lb]{\smash{\SetFigFont{6}{7.2}{rm}{\color[rgb]{0,0,0}$\ast$}%
}}}
\put(5401,-4636){\makebox(0,0)[lb]{\smash{\SetFigFont{6}{7.2}{rm}{\color[rgb]{0,0,0}$\ast$}%
}}}
\put(6001,-4636){\makebox(0,0)[lb]{\smash{\SetFigFont{6}{7.2}{rm}{\color[rgb]{0,0,0}$\ast$}%
}}}
\put(6601,-4636){\makebox(0,0)[lb]{\smash{\SetFigFont{6}{7.2}{rm}{\color[rgb]{0,0,0}$\ast$}%
}}}
\put(7201,-4636){\makebox(0,0)[lb]{\smash{\SetFigFont{6}{7.2}{rm}{\color[rgb]{0,0,0}$\ast$}%
}}}
\put(7801,-4636){\makebox(0,0)[lb]{\smash{\SetFigFont{6}{7.2}{rm}{\color[rgb]{0,0,0}$\ast$}%
}}}
\put(8401,-4636){\makebox(0,0)[lb]{\smash{\SetFigFont{6}{7.2}{rm}{\color[rgb]{0,0,0}$\ast$}%
}}}
\put(9001,-4636){\makebox(0,0)[lb]{\smash{\SetFigFont{6}{7.2}{rm}{\color[rgb]{0,0,0}$\ast$}%
}}}
\put(9601,-4636){\makebox(0,0)[lb]{\smash{\SetFigFont{6}{7.2}{rm}{\color[rgb]{0,0,0}$\ast$}%
}}}
\put(8401,-3436){\makebox(0,0)[lb]{\smash{\SetFigFont{6}{7.2}{rm}{\color[rgb]{0,0,0}$\ast$}%
}}}
\put(9001,-3436){\makebox(0,0)[lb]{\smash{\SetFigFont{6}{7.2}{rm}{\color[rgb]{0,0,0}$\ast$}%
}}}
\put(9601,-3436){\makebox(0,0)[lb]{\smash{\SetFigFont{6}{7.2}{rm}{\color[rgb]{0,0,0}$\ast$}%
}}}
\put(7801,-3436){\makebox(0,0)[lb]{\smash{\SetFigFont{6}{7.2}{rm}{\color[rgb]{0,0,0}$\ast$}%
}}}
\put(7201,-2236){\makebox(0,0)[lb]{\smash{\SetFigFont{6}{7.2}{rm}{\color[rgb]{0,0,0}$\ast$}%
}}}
\put(7801,-2236){\makebox(0,0)[lb]{\smash{\SetFigFont{6}{7.2}{rm}{\color[rgb]{0,0,0}$\ast$}%
}}}
\put(8401,-2236){\makebox(0,0)[lb]{\smash{\SetFigFont{6}{7.2}{rm}{\color[rgb]{0,0,0}$\ast$}%
}}}
\put(9001,-2236){\makebox(0,0)[lb]{\smash{\SetFigFont{6}{7.2}{rm}{\color[rgb]{0,0,0}$\ast$}%
}}}
\put(9601,-2236){\makebox(0,0)[lb]{\smash{\SetFigFont{6}{7.2}{rm}{\color[rgb]{0,0,0}$\ast$}%
}}}
\put(6601,-2236){\makebox(0,0)[lb]{\smash{\SetFigFont{6}{7.2}{rm}{\color[rgb]{0,0,0}$\ast$}%
}}}
\put(4801,-1036){\makebox(0,0)[lb]{\smash{\SetFigFont{6}{7.2}{rm}{\color[rgb]{0,0,0}$\ast$}%
}}}
\put(5401,-1036){\makebox(0,0)[lb]{\smash{\SetFigFont{6}{7.2}{rm}{\color[rgb]{0,0,0}$\ast$}%
}}}
\put(6001,-1036){\makebox(0,0)[lb]{\smash{\SetFigFont{6}{7.2}{rm}{\color[rgb]{0,0,0}$\ast$}%
}}}
\put(6601,-1036){\makebox(0,0)[lb]{\smash{\SetFigFont{6}{7.2}{rm}{\color[rgb]{0,0,0}$\ast$}%
}}}
\put(7201,-1036){\makebox(0,0)[lb]{\smash{\SetFigFont{6}{7.2}{rm}{\color[rgb]{0,0,0}$\ast$}%
}}}
\put(7801,-1036){\makebox(0,0)[lb]{\smash{\SetFigFont{6}{7.2}{rm}{\color[rgb]{0,0,0}$\ast$}%
}}}
\put(8401,-1036){\makebox(0,0)[lb]{\smash{\SetFigFont{6}{7.2}{rm}{\color[rgb]{0,0,0}$\ast$}%
}}}
\put(9001,-1036){\makebox(0,0)[lb]{\smash{\SetFigFont{6}{7.2}{rm}{\color[rgb]{0,0,0}$\ast$}%
}}}
\put(9601,-1036){\makebox(0,0)[lb]{\smash{\SetFigFont{6}{7.2}{rm}{\color[rgb]{0,0,0}$\ast$}%
}}}
\put(9601,164){\makebox(0,0)[lb]{\smash{\SetFigFont{6}{7.2}{rm}{\color[rgb]{0,0,0}$\ast$}%
}}}
\put(901,-3136){\makebox(0,0)[lb]{\smash{\SetFigFont{10}{12.0}{rm}{\color[rgb]{0,0,0}$m$}%
}}}
\put(-1799,-2986){\makebox(0,0)[lb]{\smash{\SetFigFont{10}{12.0}{rm}{\color[rgb]{0,0,0}$x$}%
}}}
\put(-749,-3886){\makebox(0,0)[lb]{\smash{\SetFigFont{10}{12.0}{rm}{\color[rgb]{0,0,0}$y$}%
}}}
\put(9601,764){\makebox(0,0)[lb]{\smash{\SetFigFont{6}{7.2}{rm}{\color[rgb]{0,0,0}$\ast$}%
}}}
\put(12451,-2611){\makebox(0,0)[lb]{\smash{\SetFigFont{10}{12.0}{rm}{\color[rgb]{0,0,0}$m'$}%
}}}
\end{picture}